\documentclass[12pt,leqno]{article}

\usepackage{amssymb}
\usepackage{amsmath}
\usepackage{graphicx}
\usepackage[usenames]{color}

\def\cald{{\mathcal{D}}}
\def\call{{\mathcal{L}}}

\def\calm{{\mathcal{M}}}
\def\calf{{\mathcal{F}}}
\def\calo{{\mathcal{O}}}
\def\calp{{\mathcal{P}}}

\def\calx{{\mathcal{X}}}

\def\rr{{\mathbb{R}}}
\def\nn{{\mathbb{N}}}
\def\9{{\infty}}
\def\lbb{{\lambda}}

\def\g{{\gamma}}

\def\vf{{\varphi}}

\def\ooo{{\Omega}}
\def\pp{{\partial}}
\def\D{{\Delta}}
\def\vp{{\varepsilon}}

\def\1{^{-1}}

\def\barr{\begin{array}}
\def\earr{\end{array}}

\def\dd{\displaystyle}
\def\bk{\bigskip }

\def\n{\noindent }

\def\vsp{\vspace*{2mm}\\ }

\def\ff{\forall }

\def\({\left(}
\def\){\right)}
\def\<{\left<}
\def\>{\right>}

\newtheorem{theorem}{Theorem}[section]
\newtheorem{proposition}[theorem]{Proposition}
\newtheorem{lemma}[theorem]{Lemma}

\newtheorem{remark}[theorem]{Remark}
\newtheorem{hypothesis}[theorem]{Hypothesis}

\title{From nonlinear Fokker-Planck  equations to~solutions of distribution dependent SDE}
\author{Viorel Barbu\thanks{Octav Mayer Institute of Mathematics of the Romanian Academy, Ia\c si, Romania} \and Michael R\"ockner\thanks{Fakult\"at f\"ur Mathematik, Universit\"at Bielefeld, D-33501 Bielefeld, Germany}\ \thanks{Academy of Mathematics and Systems Science, CAS, Beijing}}
\date{}

\begin{document}

\maketitle

\begin{abstract}   We construct weak solutions to the McKean-Vlasov SDE 
	$$dX(t)=b\(X(t),
	\dd\frac{d\call_{X(t)}}{dx}(X(t))\)dt+
\sigma\(X(t),\dd\frac{d\call_{X(t)}}{dt}(X(t))\)dW(t)$$
on $\rr^d$ for possibly degenerate diffusion matrices $\sigma$ with $X(0)$ having a given law, which has a density with respect to   Lebesgue measure, $dx$. Here $\call_{X(t)}$ denotes the law of $X(t)$.   Our approach is to first solve the corresponding nonlinear Fokker-Planck equations and then use the well known superposition principle to obtain weak solutions of the above SDE.\medskip\\
{\bf Keywords:} Fokker-Planck  equation, Kolmogorov operator, probability density, $m$-accretive operator, Wiener process.\\
{\bf Mathematics Subject Classification (2000):} 60H30, 60H10,\break 60G46, 35C99, 58J165.
\end{abstract}

\section{Introduction}\label{s1}

Recently there has been an increasing interest in distribution dependent stochastic differential equations (DDSDE for short) of type
\begin{equation}\label{1.1}
\barr{l}
dX(t)=b(t,X(t),\call_{X(t)})dt+
\sigma(t,X(t),\call_{X(t)})dW(t)\vsp
X(0)=\xi_0,\earr
\end{equation}
on $\rr^d$, where $W(t),\ t\ge0$, is an $(\calf_t)$-Brownian motion on a probability space $(\ooo,\calf,P)$ with normal filtration $(\calf_t)_{t\ge0}$. The coefficients $b,$ $\sigma$ are defined on $[0,\9)\times\rr^d\times\mathcal{P}(\rr^d)$ are $\rr^d$ and $d\times d$-matrix valued, respectively (sa\-tis\-fying conditions to be specified below). Here $\mathcal{P}(\rr^d)$ denotes the set of all probability measures on $\rr^d$. In \eqref{1.1}, $\call_{X(t)}$ denotes the law of $X(t)$ under $P$ and $\xi_0$ is an $\calf_0$-measurable $\rr^d$-valued map. Equations as in \eqref{1.1} are also referred to as McKean-Vlasov SDEs. Here we refer to the classical papers \cite{4*},  \cite{1*}, \cite{2*}, \cite{5*},  \cite{3*},  and, e.g., the more recent papers \cite{9*}, \cite{8*}, \cite{11*}, \cite{10*},
\cite{12*}, \cite{7*} and \cite{6*}.

By It\^o's formula, under quite general conditions on the coefficients, the time marginal laws $\mu_t:=\call_{X(t)},$ $t\ge0$, with $\mu_0:=\mbox{law of }\xi_0$, of the solution $X(t),\ t\ge0$, to \eqref{1.1} satisfy a {\it nonlinear} Fokker-Planck equation (FPE for short). More precisely, for all $\vf\in C^2_0(\rr^d)$ (=~all twice differentiable real-valued functions of compact support) and, for all $t\ge0$,
\begin{equation}\label{1.2}
\int_{\rr^d}\vf(x)\mu_t(dx)=\int_{\rr^d}\vf(x)\mu_0(dx)+\int^t_0\int_{\rr^d}L_{\mu_s}\vf(s,x)\mu_s(dx)ds,\end{equation}where, for $x\in\rr^d,\ t\ge0,\ a_{ij}:=(\sigma\sigma^T)_{i,j},$ $1\le i,j\le d,$
\begin{equation}\label{1.3}
\dd L_{\mu_t}\vf(t,x):=\dd\frac12\sum\limits^d_{i,j=1}\!a_{ij}(t,x,\mu_t)\,\frac{\pp^2}{\pp x_i\pp x_j}\,\vf(x)
 \dd+\!\!\sum\limits^d_{i=1}\!b_i(t,x,\mu_t)\,\frac\pp{\pp x_i}\,\vf(x),
\end{equation}
is the corresponding   Kolmogorov operator. For equations of type \eqref{1.2}, we refer the reader, e.g., to \cite{13*}.   We note that \eqref{1.2} is also shortly written as
\begin{equation}\label{1.3prim}
\pp_t\mu_t=L^*_{\mu_t}\mu_t\mbox{ with }\mu_0\mbox{ given}.\end{equation}
Hence, if one can solve \eqref{1.1}, one obtains a solution to \eqref{1.2} this way.

 In the special case where  the solutions $\mu_t$, $t\ge0$, to \eqref{1.2} have densities with respect to the Lebesgue measure $dx$, i.e., $\mu_t(dx)=u(t,x)dx,$ $t\ge0$,   \eqref{1.2} can be rewritten (in the sense of Schwartz distributions) as (cf.\cite{8a})

\begin{equation}\label{1.4}
\barr{rl}
\dd\frac\pp{\pp t}\,u(t,x)=\!\!\!
&\dd\frac12\sum^d_{i,j=1}\frac{\pp^2}{\pp x_i\pp x_j}\,[a_{ij}(t,x,u(t,\cdot)dx)u(t,x)]\vsp
&-\dd\sum^d_{i=1}\frac\pp{\pp x_i}[b_i(t,x,u(t,\cdot)dx)u(t,x)]\vsp u(0,x)=\!\!\!&u_0(x)\ \(=\dd\frac{d\mu_0}{dx}\,(x)\),\earr\end{equation}$x\in\rr^d,\ t\ge0,$ or, shortly,
\begin{equation}\label{1.5}
\pp_tu=\frac12\,\pp_i\pp_j(a_{ij}(u)u)-\pp_i(b_i(u)u),\ u(0,\cdot)=u_0.\end{equation}

In this paper, we want to go in the opposite direction, that is, we first want to solve \eqref{1.2} and, using the obtained $\mu_t$, $t\ge0,$ we shall obtain a (probabilistically) weak  solution to \eqref{1.1} with the time marginal laws of $X(t), t\ge0$, given by these $\mu_t,\ t\ge0$. It turns out that, once one has solved \eqref{1.2}, which is in general a hard task, and if one can prove some mild integrability properties for the solutions, a recent version of the so-called {\it superposition principle} by Trevisan in \cite{8} (generalizing   earlier work by Figalli \cite{14*}), in connection with a classical result by Stroock and Varadhan (see, e.g., \cite{15*}) yields the desired weak solution of \eqref{1.1} (see Section \ref{s2} below for details).

We would like to mention at this point that, by the very same result from \cite{8}, one can also easily prove that, if \eqref{1.1} has a unique solution in law, then the solution to \eqref{1.2} does not only exist as described above, but is also unique. In this paper,   however, we concentrate on existence of weak solutions to \eqref{1.1}. We shall do this in the singular case, where the coefficients in \eqref{1.1} are of {\it Nemytskii-type}, that is, we consider the following situation:
      $b_i,a_{ij}$ depend  on $\mu $ in the following way: \begin{equation}\label{1.6}
b_i(t,x,\mu ):=\bar b_i\(t,x,\dd\frac{d\mu}{dx}\,(x)\),\ \ \
a_{ij}(t,x,\mu ):=\bar a_{ij}\(t,x,\dd\frac{d\mu}{dx}\,(x)\),
\end{equation}for $t\ge0,\ x\in\rr^d,\ 1\le i,j\le d,$ where $\bar b_i,\bar a_{ij}:[0,\9)\times\rr^d\times\rr\to\rr$, are measurable functions. Then, under the conditions on $\bar b_i$ and $\bar a_{ij}$, $1\le i,j\le d$, specified in Section {\rm\ref{s3}}, we shall construct solutions $(\mu_t)_{t\ge0}$ to \eqref{1.1} which are absolutely continuous with respect to the Lebesgue $dx$, i.e., $\mu_t(dx)=u(t,x)dx,$ $t\ge0$. So, as indicated above, by the superposition principle, we obtain weak solutions to DDSDEs of type
\begin{equation}\label{1.7}
\barr{rl}
dX(t)=\!\!\!
&\bar b\(t,X(t),\dd\frac{d\call_{X(t)}}{dx}\,(X(t))\)dt\vsp
&+\bar\sigma\(t,X(t),\dd\frac{d\call_{X(t)}}{dx}\,(X(t))\)dW(t),\vsp X(0)=\!\!\!&\xi_0,\earr\end{equation}with $(\bar\sigma\bar\sigma^T)_{i,j}=\bar a_{ij}$.

In particular, we obtain a probabilistic representation of the solution $\mu_t$, $t\ge0$, of the nonlinear FPE \eqref{1.2} (or \eqref{1.4}) as the time marginal laws of a stochastic process, namely the solution of the DDSDE \eqref{1.7}.

 We would like to emphasize that the coefficients as in \eqref{1.7}, which we consider below, have no continuity properties with respect to their dependence on the law $\call_{X(t)}$ of $X(t)$, such as those imposed in the existing literature on the subject. Nevertheless, such {\it Nemytskii-type}-dependence is very natural and, of course, independent of the $dx$-version of the Lebesgue density of $\call_{X(t)}$ we choose in \eqref{1.7}, since we are looking only for solutions of \eqref{1.7} in the class with $\call_{X(t)}$ being absolutely continuous with respect to $dx$. Precise conditions on the coefficients $\bar b_i,\bar a_{ij}$ are formulated in Section \ref{s3} (there, for simplicity, denoted by $b_i,a_{ij}$). Our main existence results for   solutions of the nonlinear FPE \eqref{1.2} are Theo\-rems \ref{t1} and \ref{t2} below. Our main result  on solutions to \eqref{1.1} (more precisely, \eqref{1.7}) is Theorem  \ref{t4.1}.  Subsequently, in Remark \ref{r4.1}  we discuss connections with previous related, but much more special,  results from \cite{3}--\cite{5}.  A class of cases where we also have uniqueness in law results for solutions to \eqref{1.7} is described in Remark \ref{r4.2}.

\medskip\noindent{\bf Notations.}~Given an open subset $\calo\subset\rr^d$, we denote by $L^p(\calo)$, \mbox{$1\le p\le\9$,}\break the standard Lebesgue $p$-integrable functions on $\calo$, and by $H^1(\calo)$, the Sobolev space $\{u\in L^2(\calo);\nabla u\in L^2(\calo)\}.$

We set $H^1_0(\calo)=\{u\in H^1(\calo);$ $ u=0\mbox{ on }\pp\calo\}$ and denote by $H\1(\calo)$ the dual space of $H^1_0(\calo)$.
By $C^\9_0(\calo)$ we denote the space of infinitely differentiable functions with compact support in $\calo$. We set $H^1=H^1(\rr^d)$, $H\1=H\1(\rr^d)$ and denote by $H^1_{\rm loc}$ the corresponding local space.

We   also set $L^p=L^p(\rr^d)$ with the norm denoted $|\cdot|_p$ and $L^p_{\rm loc}=L^p_{\rm loc}(\rr^d)$, $1\le p\le\9.$
 By $\cald'(\rr^d)$ and $\cald'((0,\9)\times\rr^d)$, we denote the space of distributions on $\rr^d$ and $(0,\9)\times\rr^d$, respectively.
 
 We shall denote either by $\frac{\partial u}{\partial x_j}$ or by $u_{x_j},D_j$, $u$  the partial derivative of the function $u=u(x_1,...,x_d)$ with respect to $x_j,$ $1\le j\le d.$ By $D^2_{ij},$ $i,j=1,...,d,$ we shall denote the second order derivatives $\frac{\partial^2u}{\partial x_i\partial x_j}$.

We denote by $C(\rr^d\times\rr)$ and $C(\rr^d)$ the space of continuous functions on $\rr^d\times\rr$ and $\rr$, respectively, and by $C_b(\rr^d\times\rr)$ and $C_b(\rr)$  the corresponding subspaces of  continuous and bounded functions.

By $C^1(\rr^d\times\rr)$ and $C^1(\rr)$ we denote the spaces of continuously differentiable functions on $\rr^d\times\rr$ and $\rr$, respectively. 
Finally, $C^1_b$ is the space of bounded continuously differentiable functions with bounded derivatives.

If $\calx$ is a real Banach space and $0<T<\9$, we denote by $L^p(0,T;\calx)$ the space of Bochner $p$-integrable functions $u:(0,T)\to\calx$ and by $C([0,T];\calx)$ the space of of $\calx$-valued continuous functions on $[0,T]$.

\section{From nonlinear FPEs to DDSDEs:\\ general scheme}\label{s2}
\setcounter{equation}{0}

Let $a_{ij},b_i:[0,\9)\times\rr^d\times\calp(\rr^d)\to\rr,$ $1\le i,j\le d$, be measurable.

\begin{hypothesis}\label{h2.1} There exists a solution $(\mu_t)_{t\ge0}$ to \eqref{1.2} such that
	\begin{itemize}
		\item[\rm(i)] $\mu_t\in\calp(\rr^d)$ for all $t\ge0$.
		\item[\rm(ii)] $(t,x)\mapsto a_{ij}(t,x,\mu_t)$ and $(t,x)\mapsto b_i(t,x,\mu_t)$ are measurable and
		 $$\int^T_0\int_{\rr^d}[|a_{ij}(t,x,\mu_t)|+|b_i(t,x,\mu_t)|]\mu_t(dx)dt<\9,$$for all $T\in(0,\9).$
		\item[\rm(iii)] $[0,\9)\ni t\mapsto\mu_t$ is weakly continuous.
	\end{itemize}
\end{hypothesis}
Under Hypothesis \ref{h2.1}, we can apply the superposition principle (see Theorem 2.5 in \cite{8})  for {\it linear} FPEs applied to the ({\it linear}) Kolmogorov operator
\begin{equation}\label{e2.1z}
L_{\mu_t}:=\frac12\sum^d_{i,j=1}a_{ij}(t,x,\mu_t)\frac{\pp^2}{\pp x_i\pp x_j}+\sum^d_{i=1}b_i(t,x,\mu_t)\frac\pp{\pp x_i},
\end{equation}with   $(\mu_t)_{t\ge0}$ from Hypothesis \ref{h2.1} fixed.

More precisely, by Theorem 2.5 in \cite{8}, there exists a probability measure $P$ on $C([0,T];\rr^d)$ equipped with its Borel $\sigma$-algebra and its natural normal filtration obtained by the evaluation maps $\pi_t$, $t\in[0,T]$, defined by
$$\pi_t(w):=w(t),\ w\in C([0,T],\rr^d),$$solving the martingale problem  (see \cite{8}, Definition 2.4) for the time-dependent ({\it linear}) Kolmogorov ope\-ra\-tor  $\frac\pp{\pp t}+L_{\mu_t}$ (with $(\mu_t)_{t\ge0}$ as above fixed) with time marginal  laws $$P\circ\pi^{-1}_t=\mu_t,\ t\ge0.$$ Then, a standarad result (see Theorem 4.5.2 in \cite{15*}) implies that there exists a $d$-dimensional $(\calf_t)$-Brownian motion $W(t)$, $t\ge0$, on a stochastic basis $(\ooo,\calf,(\calf_t)_{t\ge0},Q)$ and a continuous $(\calf_t)$-progressively measurable map $X:[0,\9)\times\ooo\to\rr^d$ satisfying the following (DD)SDE
\begin{equation}\label{e2.2z}
dX(t)=b(t,X(t),\mu_t)dt+\sigma(t,X(t),\mu_t)dW(t),\end{equation}
with the law
$$Q\circ X\1=P,$$
where $\sigma=((a_{ij})_{1\le i,j\le d})^{\frac12}.$ In particular, we have, for the marginal laws,
\begin{equation}\label{e2.3z}
\call_{X(t)}:=Q\circ X(t)\1=\mu_t,\ t\ge0.\end{equation}

\begin{remark}\label{r2.1}\rm Because of \eqref{e2.3z}, the process $X(t),\ t\ge0$, is also called a {\it probabilistic representation} of the solution $(\mu_t)_{t\ge0}$ for the nonlinear FPE \eqref{1.2}.\end{remark}

\begin{remark}\label{r2.2} \rm It is much harder to prove that the solution to SDE  (\ref{e2.2z}) for fixed $(\mu_t)_{t\ge0}$ is unique in law, provided its initial distribution is $\mu_0$, which would, of course, be very desirable. For this, one has to prove the uniqueness of the solutions to the {\it linear} Fokker-Planck equation
	$$\pp_t\nu_t=L^*_{\mu_t}\nu_t,\ \nu_0=\mu_0$$for all initial condition of the type $\mu_0=\delta_x$, $x\in\rr^d$  (see \cite{15*}, Theorem 6.2.3).   For a large class of initial conditions $\mu_0$, this was achieved in certain cases where $d=1$ (see \cite{4}, \cite{5}, \cite{13*a}). \end{remark}

\paragraph{\it Conclusion.} To weakly solve  DDSDE \eqref{1.1}, we have to solve the corres\-pon\-ding nonlinear FPE \eqref{1.2} (hard!) and then check Hypothesis \ref{h2.1} above.

\section{Existence of solutions to the nonlinear~FPEs}
\label{s3}
\setcounter{equation}{0}

Consider the following time-independent special case of \eqref{1.4} with  Nemytskii-type dependence of the coefficients on $u(t,x)dx,$ $t\ge0$, i.e., the nonlinear Fokker-Planck equation
\begin{equation}\label{e1.1}
\barr{l}
\dd\frac{\pp u}{\pp t}-\sum^d_{i,j=1}D^2_{ij}(a_{ij}(x,u)u)+{\rm div}(b(x,u)u)=0\mbox{ in }\cald'((0,\9)\times\rr^d),\vsp
u(0,x)=u_0(x),\ x\in\rr^d,\earr\end{equation}where $b(x,u)=\{b_i(x,u)\}^d_{i=1}$.\medskip

We shall study this equation under two different sets of hypotheses spe\-ci\-fied in the following.
\begin{itemize}
	\item[(H1)] $a_{ij}\in C^2(\rr^d\times\rr)\cap C_b(\rr^d\times\rr^d),\
	(a_{ij})_{x}\in C_b(\rr^d\times\rr;\rr^d),\   \mbox{$\ff a_{ij}=a_{ji},$}  $ $ i,j=1,...,d.$
	\item[(H2)] $\sum\limits^d_{i,j=1}(a_{ij}(x,u)+(a_{ij}(x,u))_uu)\xi_i\xi_j\ge\g|\xi|^2,\ \ff\xi\in\rr^d,\ x\in\rr^d,\  u\in\rr,$ where $\g>0.$
	\item[(H3)] $b_i\in C_b(\rr^d\times\rr)\cap C^1(\rr^d\times\rr),$ $b_i(x,0)\equiv0,$ $\ff x\in\rr^d,\ i=1,...,d.$
	\item[(H1)$'$] $a_{ij}(x,u)\equiv a_{ij}(u),\  a_{ij}\in C^2(\rr)\cap C_b(\rr),$ $a_{ij}=a_{ji},$ $\ff i,j=1,...,d.$
	\item[(H2)$'$] $\sum\limits^d_{i,j=1}(a_{ij}(u)+u(a_{ij}(u))_u)\xi_i\xi_j\ge0,\ \ff\xi\in\rr^d,\ u\in\rr.$
	\item[(H3)$'$] $b_i\in C_b(\rr)\cap C^1(\rr^d),\ b_i(0)=0,\ i=1,...,d.$
\end{itemize}Here $(a_{ij}(x,u))_u=\frac\pp{\pp u}\,a_{ij}(x,u),\ \ff u\in\rr,$ and  $(a_{ij})_x(x,u)=(\nabla_x a_{ij})(x,u)$, $x=\{x_i\}^d_{i=1}$.  
The first set of hypotheses, that is (H1)--(H3), allows for nonlinear non\-de\-ge\-nerate FPEs   with $x$-dependent coefficients, while the se\-cond set (H1)$'$--(H3)$'$ allows for degenerate nonlinear FPEs, however, with $x$-independent coefficients.

Nonlinear FPEs of the form \eqref{e1.1} describe in the mean field theory the dynamics of a set of interacting particles or many body systems. The function $u=u(t,x)$ is associated with the probability to find a certain subsystem or particle at time $t$ in the state $x$. Equation \eqref{e1.1} arises also as a closed loop system corresponding to a velocity field system
$$\frac{\pp v}{\pp t}=F(x,u)v=\sum^d_{i,j=1}D^2_{ij}(a_{ij}(x,u)v)-{\rm div}(b(x,u)v)$$
with coefficients depending on the probability density $u$. If $v=u$, one may view this system   as a statistical feedback (see \cite{8a}).

The first part of this section is concerned with the existence of a weak (mild) solution to equation \eqref{e1.1}  in the space $L^1(\rr^d)$. This result is obtained via the Crandall and Liggett existence theorem for the nonlinear Cauchy problem
\begin{equation}
\barr{l}
\dd\frac{du}{dt}\,(t)+Au(t)=0,\ t\ge0,\vsp
u(0)=u_0,\earr\label{e1.2}
\end{equation}
in a Banach space $\calx.$

An operator $A:D(A)\subset\calx\to\calx$ (possibly multivalued) is said to be $m$-accretive if, for each $\lbb>0$, the range $R(I+\lbb A)$ of the operator $I+\lbb A$ is all of $\calx$ and
\begin{equation}\label{e1.3}
\|(I+\lbb A)\1u-(I+\lbb A)\1v\|_\calx\le\|u-v\|_\calx,\ \ff u,v\in\calx,\ \lbb>0.\end{equation}
The continuous function $u:[0,\9)\to\calx$ is said to be a mild solution to \eqref{e1.2} if, for each $0<T<\9,$
\begin{eqnarray}
&&u(t)=\lim\limits_{h\to0} u_h(t)\mbox{ strongly in $\calx$, uniformly in $t\in[0,T]$}\hspace*{18mm}\label{e1.4}\\
&&\mbox{where $u_h:[0,T]\to\calx$ is defined by}\nonumber\\
&&u_h(t)=u^i_h,\ t\in[ih,(i+1)h),\ i=0,1,...,N=\left[\frac Th\right].\label{e1.5}\\
&&u^i_h+hAu^i_h=u^{i-1}_h,\ i=1,...,N;\ u^0_h=u_0.\label{e1.6}
\end{eqnarray}
By the Crandall and Liggett theorem (see, e.g.,  \cite{1}, p.~99), if $A$ is $m$-accretive, then for each $u_0\in\overline{D(A)}$ (the closure of $D(A)$ in $\calx)$   there is a unique mild solution $u\in C([0,\9);\calx)$ to \eqref{e1.2}.  Moreover, the map $u_0\to u(t)$ is a continuous semigroup of contractions on $\overline{D(A)}$ equipped with $\|\cdot\|_\calx$.

The first main existence result of this section, Theorem \ref{t1}, is obtained by \mbox{writing} equation \eqref{e1.1} in the form \eqref{e1.2} with a suitable $m$-accretive operator $A$ in the space $\calx=L^1(\rr^d)$.

It should be said that the space $L^1(\rr^d)$ is not only appropriate to re\-pre\-sent equation \eqref{e1.1} in the form \eqref{e1.2}, but it is the unique $L^p(\rr^d)$-space in which the operator defined by equation \eqref{e1.1} is $m$-accretive, that is, which gives the parabolic character of this equation. Only in the particular case of porous media equations (i.e., \eqref{e1.1} with $b\equiv0$), an alternative is the Sobolev space $H\1(\rr^d)$, but this does not work for the more general case \eqref{e1.1}. On the other hand, taking into account the significance of the solution $u$ as probability density, the space $L^1$ is very convenient for the treatment of equation \eqref{1.1}.

Our work \cite{3} contains the following     special case of \eqref{e1.1}:
\begin{equation}\label{e1.7}
 \frac{\pp u}{\pp t}-\D\beta(u)+{\rm div}(b(u)u)=0\mbox{ in }(0,T)\times\rr^d,
\end{equation}where $\beta:\rr\to2^\rr$ is a maximal monotone (multivalued) function with $\sup\{|s|:s\in\beta(r)\}\le C|r|^m,\ r\in\rr,$ for some $C,m\in[0,\9)$. (See also~\cite{2}.)  In the special case $b\equiv0$ and $d=1$, related results were obtained in \cite{4}, \cite{5}. However,  the present case is much more difficult and the arguments of \cite{3} are not applicable here.

\subsection{Existence for  FPEs in the nondegenerate,\\ $x$-dependent case}

Define in the space $\calx=L^1$ the operator $A:D(A)\subset L^1\to L^1$,
\begin{eqnarray}
Au&=&-\sum^d_{i,j=1}D^2_{ij}(a_{ij}(x,u)u)+{\rm div}(b(x,u)u),\ \ff u\in D(A),\label{e2.1}\\
D(A)&=&\left\{u\in L^1;-\sum^d_{i,j=1}D^2_{ij}(a_{ij}(x,u)u)+{\rm div}(b(x,u)u)\in L^1\right\},\label{e2.2}
\end{eqnarray}where $D^2_{ij}$ and div are taken in the sense of Schwartz distributions on $\rr^d$, i.e., in $\cald'(\rr^d)$. 
We note that, since by (H1), (H3), $a_{ij}(x,u)u$, $b_i(x,u)u\in L^1$, $\ff i,j=1,...,d$,  $\ff u\in L^1$, $Au$ is well defined in $\cald'(\rr^d)$. Moreover, since $C^\infty_0(\rr^d)\subset D(A)$, it follows that $D(A)$ is dense in $L^1$.

Since we are going to represent equation \eqref{e1.1} as \eqref{e1.2} with $A$ defined by \eqref{e2.1}--\eqref{e2.2}, we must prove that $A$ is $m$-accretive, that is, $R(I+\lbb A)=L^1$ and \eqref{e1.3} holds in $\calx=L^1$ for all $\lbb>0.$ For this purpose, we shall prove the following result.

\begin{proposition}\label{p1} Let {\rm(H1)--(H3)} hold. Then, for each $f\in L^1$ and $\lbb>0$, the equation
\begin{equation}\label{e2.3}
u-\lbb\sum^d_{i,j=1}D^2_{ij}(a_{ij}(x,u)u)+\lbb\ {\rm div}(b(x,u)u)=f\mbox{ in }\cald'(\rr^d)\end{equation}has a unique solution $u=u(\lbb,f)\in D(A).$ 

 Moreover,  we have, for all $\lbb>0,$ 
\begin{eqnarray}
&&|u(\lbb,f_1)-u(\lbb,f_2)|_1\le|f_1-f_2|_1,\ \ff f_1,f_2\in L^1,\label{e2.5}\\
&&(I+\lbb A)\1f\ge0,\mbox{ a.e. in }\rr^d,\mbox{ if }f\in L^1,\ f\ge0,\ \mbox{ a.e. in }\rr^d,\label{e2.7}\\
&&\int_{\rr^d}(I+\lbb A)\1f(x)dx=\int_{\rr^d}f(x)dx,\ \ff f\in L^1,\ \lbb>0.\label{e2.8}
\end{eqnarray}
\end{proposition}

\noindent{\bf Proof.}
In the following, we shall simply write
$$a_{ij}(u)=a_{ij}(x,u),\ x\in\rr^d,\ u\in\rr,\ i,j=1,...,d.$$
We set
$$\barr{c}
a^*_{ij}(u)\equiv a_{ij}(x,u)u,\ x\in\rr^d,\ u\in\rr,\ \ff i,j=1,...,d,\vsp 
b(x,u)=\{b_i(x,u)\}^d_{i=1},\ b^*(x,u)=b(x,u)u,\ x\in\rr^d,\ u\in\rr.\earr$$We note that, by   (H2), we have
\begin{equation}
\sum^d_{i,j=1}(a^*_{ij})_u(x,u)\xi_i\xi_j\ge\g|\xi|^2,\ \ff\xi\in\rr^d,\ x\in\rr^d,\ u\in\rr,\label{e2.11}
\end{equation}where $\gamma>0.$ 
We shall first prove Proposition \ref{p1} under the additional hypotheses
\begin{itemize}
	\item[\rm(K)] $(a^*_{ij})_u\in C_b(\rr^d\times\rr), b_i\in C^1_b(\rr^d\times\rr)$, and\end{itemize}
	\begin{equation}\label{e3.16az}
	\barr{r}
	|(a^*_{ij})_u(x,u)-(a^*_{ij})_u(x,\bar u)|  +
	|\nabla_x(a^*_{ij})(x,u)-\nabla_x(a^*_{ij})(x,\bar u)|_d\vsp
	\le C|u-\bar u|,\  \ff u,\bar u\in\rr,\ x\in\rr^d,\earr\end{equation}
	\begin{itemize}\item[\ ]for $i,j=1,...,d.$\end{itemize}

We rewrite \eqref{e2.3} as
$$u-\lbb\sum^d_{i,j=1}D^2_{ij}(a^*_{ij}(u))+\lbb\ {\rm div}(b^*(x,u))=f\mbox{ in }\cald'(\rr^d).\leqno\eqref{e2.3}'$$
Equivalently, if $Du\in L^1_{\rm loc},$ then
$$\begin{array}{r}
u-\lbb\sum\limits^d_{i,j=1}D_{i}
((a^*_{ij})_u(u)D_ju+(a^*_{ij})_{x_j}(x,u)u)+\lbb\ {\rm div}(b^*(x,u))=f\\\mbox{ in }\cald'(\rr^d).\end{array}\leqno\eqref{e2.3}''$$
We also set
$$\barr{rcl}
b_\9&=&\sup\{|b_i(x,u)|;\ (x,u)\in\rr^d\times\rr,\ i=1,...,d\},\vsp
c_\9&=&\sup\{|(a_{ij})_{x_j}(x,u)|;\ (x,u)\in\rr^d\times\rr,\ i,j=1,...,d\}.\earr$$
(By virtue of (K), the   formulation  $\eqref{e2.3}''$ of $\eqref{e2.3}'$ makes sense only if $D_ju\in L^1_{\rm loc}.)$

For each $N>0$, we set $B_N=\{\xi\in\rr^d;\ |\xi|< N\}$. We have

\begin{lemma}\label{l1} Let $f\in L^2$ and $0<\lbb<\lbb_0=\g(b^2_\9+c^2_\9)\1.$   Then, for each $N$  there is at least one solution $u_N\in H^1_0(B_N)$ to the equation
\begin{equation}\label{e2.12}
\barr{l}
u-\lbb\dd\sum^d_{i,j=1}D^2_{ij}(a^*_{ij}(u))+\lbb\ {\rm div}(b^*(x,u))=f\mbox{ in }B_N,\vsp
u=0\mbox{ on }\pp B_N,\earr\end{equation}	which satisfies the estimate
\begin{equation}\label{e2.13}
\|u_N\|^2_{L^2(B_N)}+\lbb\g\|\nabla u_N\|^2_{L^2(B_N)}\le C\|f\|^2_{L^2(B_N)},\end{equation}where $C$ is independent of $N$ and $\lbb.$\end{lemma}

\noindent{\bf Proof.} For $\rho>0$, we set $\calm_\rho=\{v\in L^2(B_N);\ \|v\|_{L^2(B_N)}\le\rho\}$ and consider the operator $F:\calm_\rho\to L^2(B_N)$ defined by $F(v)=u\in H^1_0(B_N)$, where $u$ is the solution to the linear elliptic problem
\begin{equation}\label{e2.14}
\barr{l}
u{-}\lbb\dd\sum^d_{i,j=1}D_i((a^*_{ij})_v(x,v)D_ju
{+}(a_{ij})_{x_j}(x,v)u)
{+}\lbb\ {\rm div}(b(x,v)u)=f\vspace*{-3mm}\\\hfill\mbox{ in }\cald'(B_N),\\
u=0\mbox{ on }\pp B_N.\earr\end{equation}
By \eqref{e2.11} and (H2), it follows via the Lax-Milgram lemma that, for each \mbox{$v\in\calm_\rho$} and $\lbb\in(0,\lbb_0)$, problem \eqref{e2.14} has a unique solution $u=F(v)$.   Moreover, by \eqref{e2.14} and (H1), we see that
\begin{equation} 
\barr{lcl}
\|u\|^2_{L^2(B_N)}+\g\lbb\|\nabla u\|^2_{L^2(B_N)}\vsp
 \qquad\le  \lbb b_\9\|\nabla u\|_{L^2(B_N)}
\|u\|_{L^2(B_N)}{+}c_\9\lbb\|u\|_{L^2(B_N)}\|\nabla u\|_{L^2(B_N)} \vsp 
 \qquad+ \|f\|_{L^2(B_N)}\|u\|_{L^2(B_N)}\label{e2.15}\vsp
 \qquad \le  \lbb b_\9\|\nabla u\|_{L^2(B_N)}\rho+\rho\|f\|_{L^2(B_N)}
+c_\9\lbb\rho\|\nabla u\|_{L^2(B_N)}.\earr\hspace*{-5mm}\end{equation}Hence, for $\lbb\in(0,\lbb_0)$ and $\rho$ suitably chosen, independent of $N$, \mbox{$F(\calm_\rho){\subset}\calm_\rho.$}

Indeed, if $v_n\to v$ in $L^2(B_N)$ and $u_n=Fv_n$, we have $b(x,v_n)\to b(x,v),$
$$(a^*_{ij})_v(x,v_n)\to (a^*_{ij})_v(x,v),\
(a_{ij})_{x_j}(x,v_n)\to (a_{ij})_{x_j}(x,v)$$
strongly in $L^2(B_N)$. Along a subsequence we have, by \eqref{e2.15},
$$u_n\to u\mbox{ weakly in }H^1(B_N),\mbox {strongly in }L^2(B_N).$$Now, letting $n\to\9$ in equation \eqref{e2.14}, where $v=v_n$ and $u=u_n$, that is,
$$\barr{r}
\dd\int_{B_N}\Big(u_n\psi+\lbb\sum\limits^d_{i,j=1}(a^*_{ij})_v(x,v_n)
-D_ju_n+(a_{ij})_{x_j}(x,v_n),u_n\Big)D_i\psi dx\vsp-\lbb\dd\int_{B_N}u_nb^*(x,v_n)\cdot\nabla\psi dx
=\dd\int_{B_N}f\psi dx,\ \ff\psi\in C^\9_0(B_N),\earr$$we see that $u=Fv$ and, therefore, $F$ is continuous on $L^2(B_N)$. 

Moreover, since the Sobolev space $H^1(B_N)$ is compactly embedded in $L^2(B_N)$, by \eqref{e2.15}  we see that $F(\calm_\rho)$ is relatively compact in $L^2(B_N)$. Then, by the Schauder theorem, $F$ has a fixed point  $u_N\in\calm{_\rho}$  which, clearly, is a solution to \eqref{e2.12}. Also, by \eqref{e2.15}, it follows that estimate \eqref{e2.13} holds.

\begin{lemma}\label{l2} Let $f\in L^2(\rr^d)$ and $\lbb<\lbb_0$. Then equation \eqref{e2.3} has at least one solution $u\in H^1(\rr^d)$ which satisfies the estimate
\begin{equation}\label{e2.16}
|u|^2_2+\g\lbb |\nabla u|^2_2\le C (|f|^2_2+1).
\end{equation}
\end{lemma}

\noindent{\bf Proof.} Consider a sequence $\{N\}\to\9$ and $u_N\in H^1_0(B_N)$ a solution to \eqref{e2.12} given by Lemma \ref{l1}. By \eqref{e2.13}, we have
$$\|u_N\|_{H^1_0(B_N)}\le C,\ \ff N,$$and so, on a subsequence, again denoted $\{N\}$, we have
\begin{equation}\label{e2.17}
u_N\to u\mbox{ weakly in } H^1(\rr^d), \mbox{ strongly in }L^2_{\rm loc}(\rr^d).\end{equation}
	Then, letting $N\to\9$ in the equation	 $$u_N{-}\lbb\!\sum^d_{i,j=1}\!\!D_i
((a^*_{ij})_uu_N)D_ju_N{+}((a_{ij})_{x_j}(x,u_N)u_N)
{+}\lbb\, {\rm div}(b(x,u_N)u_N){=}f\mbox{ in }B_N,$$or, more precisely, in its weak form\newpage
$$\barr{ll}
\dd\int_{\rr^d}u_N\psi\,dx\!\!\!
&+\lbb\dd\sum^d_{i,j=1}\int_{\rr^d}(a^*_{ij}(u_N)D_ju_N
+(a_{ij})_{x_j}(x,u_N)u_N)
D_i \psi\,dx\\
&-\lbb\dd\sum^d_{i=1}\int_{\rr^d}b(x,u_N)u_N\cdot\nabla\psi\,dx=0,\ \ff\psi\in C^\9_0(\rr^d),\earr$$
we infer by (H1), (H3) and \eqref{e2.17} that $u\in H^1(\rr^d)$ is a solution to \eqref{e2.3}. Also, estimate \eqref{e2.16} follows by \eqref{e2.13}. 
This completes the proof of Lemma~\ref{l2}.\bk

Now, we come back to the proof of Proposition \ref{p1}. We prove first that, for each $f\in L^2\cap L^1$ and $\lbb\in(0,\lbb_0)$, the solution $u=u(\lbb,f)\in H^1$  to equation \eqref{e2.3} is unique and  we have
\begin{equation}\label{e2.19}
|u(\lbb,f_1)-u(\lbb,f_2)|_1\le|f_1-f_2|_1,\ \ff f_1,f_2\in L^2\cap L^1.\end{equation}Here is the argument. We set $u_i=u(\lbb,f_i)$, $i=1,2,$ and $f=f_1-f_2$, $u=u_1-u_2$. Then, we have
\begin{equation}\label{e2.20}
\barr{l}\dd
u-\lbb\sum^d_{i,j=1}D^2_{ij}(a^*_{ij}(x,u_1)-a^*_{ij}(x,u_2))\\
\qquad\qquad\ +\lbb\ {\rm div}(b^*(x,u_1)-b^*(x,u_2))=f\mbox{ in }\cald'(\rr^d).\earr\end{equation}
More precisely, since $u_i\in H^1(\rr^d)$, equation \eqref{e2.20} is taken in its weak form
\begin{equation}\label{e3.24prim}
\barr{ll}
\dd\int_{\rr^d}\Big(u\psi+\lbb\sum^d_{i,j=1}
D_i(a^*_{ij}(x,u_1)
-a^*_{ij}(x,u_2))D_j\psi\\
\quad-\lbb\dd(b^*(x,u_1)-b^*(x,u_2))
\cdot\nabla\psi\Big)dx =\dd\int_{\rr^d}f\psi\,dx,\  \ff\psi\in H^1(\rr^d),\earr\end{equation}
In order to fix the idea of the proof, we invoke first a heuristic argument. Namely, if we multiply
\eqref{e2.20} by $\eta\in L^\9(\calo)$, $\eta(x)\in {\rm sign}(u(x))$, a.e. $x\in\rr^d$, and take into account that, by the monotonicity of the functions $a^*_{ij}$,
$$\eta(x)\in{\rm sign}(a^*_{ij}(x,u_1(x))-a^*_{ij}(x,u_2(x))),\mbox{ a.e. }x\in\rr^d,$$we get  
$$\barr{c}
\dd|u|_1+\lbb\int_{\rr^d}\sum^d_{i,j=1}D_i(a^*_{ij}(x,u_1(x))-a^*_{ij}(x,u_2(x)))D_j\eta(x)dx\\\dd
+\lbb\int_{\rr^d}{\rm div}(b^*(x,u_1)-b^*(x,u_2))\eta\,dx=\int_{\rr^d}f\eta\,dx.\earr$$
\newpage
Taking into account that, by the monotonicity of $u\to a^*_{ij}(x,u)$, we have  (formally)
$$
D_i(a^*_{ij}(x,u_1(x))-a^*_{ij}(x,u_2(x)))D_j\eta(x)\ge0\mbox{ in }\rr^d,$$
while
$$((a^*_{ij})_{x_i}(x,u_1)-(a^*_{ij})_{x_i}(x,u_2))D_j\eta(x)=0,$$
$$ \int_{\rr^d}{\rm div}(b^*(x,u_1)-b^*(x,u_2))\eta\,dx{=}
\int_{[|u|=0]}(b^*(x,u_1)-b^*(x,u_2))\cdot\nabla\eta\,dx=0,$$we get \eqref{e2.19}.
This formal argument can be made rigorous by using a smooth approximation $\calx_\delta$ of the signum graph. Namely, let $\calx_\delta\in{\rm Lip}(\rr)$ be the function
$$\calx_\delta(r)=\left\{\barr{rl}
1&\mbox{ for }r\ge\delta,\vsp\dd\frac r\delta&\mbox{ for }|r|<\delta,\vsp-1&\mbox{ for }r<-\delta,\earr\right.$$
where $\delta>0$. We note that, since $u\in L^2$, it follows $\calx_\delta(u)\in L^2$ and 
\begin{equation*}\label{e3.24az}
\barr{ll}
G_\vp=\!\!\!&-\calx_\delta(u)   \Big(\dd\sum^d_{i,j=1}D_iD_j(a^*_{ij}(x,u_1){-}a^*_{ij}(x,u_2))\\
&+{\rm div}(b^*(x,u_1)-b^*(x,u_2))\Big)\!\in\! L^1,\earr\end{equation*}
 and, therefore,
 $$\barr{r}
 \dd\int_{\rr^d}G_\vp dx = \int_{\rr^d} 
 \Big(\sum^d_{i,j=1}\!D_j(a^*_{ij}(x,u_1)-a^*_{ij}(x,u_2))D_i\calx_\delta(u)\\
  +(b^*(x,u_1)-b^*(x,u_2))\cdot\nabla
 \calx_\delta(u)\Big)dx.\earr$$
 Since $b^*(x,u_i)\in L^2$ and since, by (K), it follows that $a^*_{ij}(x,u_i)\in H^1(\rr^d)$, $i=1,2,$  taking in \eqref{e3.24prim} $\psi=\mathcal{X}_\delta(u)$  yields   
 
\begin{equation}\label{e2.21}
\barr{l}
\dd\int_{[|u(x)|\ge\delta]}|u(x)|dx+\frac1\delta\int_{[|u(x)|\le\delta]}|u(x)|dx\vsp
\quad+\lbb\dd\sum^d_{i,j=1}\int_{\rr^d}D_j(a^*_{ij}(x,u_1)-a^*_{ij}(x,u_2))D_i(\calx_\delta(u))dx\vsp
\quad=\lbb\dd\int_{\rr^d}
(b^*(x,u_1)-b^*(x,u_2))\cdot\nabla(\calx_\delta(u))dx+\int_{\rr^d}f\calx_\delta(u)dx.\earr
\end{equation}
We set
$$\barr{ll}
I^1_{\delta}\!\!\!&=\dd\int_{\rr^d}(b^*(x,u_1)-b^*(x,u_2))\cdot\nabla(\calx_\delta(u))dx\vsp
&=\dd\int_{\rr^d}(b^*(x,u_1)-b^*(x,u_2))
\cdot\nabla u \calx'_\delta(u)dx\vsp
&=\dd\frac1\delta
\int_{[|u|\le\delta]}(b^*(x,u_1)-b^*(x,u_2))\cdot\nabla u\,dx.\earr$$
Since, by hypothesis (K), $|b^*(x,u_1)-b^*(x,u_2)|\le C|u|(|u_1|+|u_2|)$, a.e. on $\rr^d$, and $u_i\in L^2$, it follows  that
$$\lim_{\delta\to0}\frac1\delta\int_{[|u|\le\delta]} 
|(b^*(x,u_1)-b^*(x,u_2))\cdot\nabla u|dx\le C
\lim_{\delta\to0}\(\int_{[|u|\le\delta]} |\nabla u|^2dx\)^{\frac12}=0,$$because $u\in H^1(\rr^d)$ and $\nabla u=0$ on $[x;u(x)=0]$.   This yields
\begin{equation}\label{e2.22}
 \lim_{\delta\to0}I^1_{\delta}=0.
\end{equation}On the other hand, taking into account that  $u_i,a^*_{ij}(u_i)\in H^1(\rr^d)$, for $i=1,2,$ we have
\begin{equation}\label{e2.22a}
\barr{ll}
I^2_{\delta}\!\!\!&=\dd\int_{\rr^d}
\sum^d_{i,j=1}D_j(a^*_{ij}(x,u_1)-a^*_{ij}(x,u_2))
D_i(\calx_\delta(u))dx\vsp
&=\dd\frac1\delta
\dd\int_{E_\delta}\sum^d_{i,j=1}
((a^*_{ij})_u(x,u_1) 
D_ju_1-
(a^*_{ij})_u(u_2) D_ju_2\vsp
&\quad+(a^*_{ij})_{x_j}(x,u_1) -
(a^*_{ij})_{x_j}(x,u_1))D_iu\,dx\vsp
&=\dd\frac1\delta\dd\int_{E_\delta}\sum^d_{i,j=1}
(a^*_{ij})_u(x,u_1) D_juD_iu\,dx\vsp
&\quad+\dd\frac1\delta\dd\int_{E_\delta}\sum^d_{i,j=1}
((a^*_{ij})_u(x,u_1)-(a^*_{ij})_u(x,u_2))
D_ju_2 D_iu\,dx\vsp
&\quad+\dd\frac1\delta\dd\int_{E_\delta}\sum^d_{i,j=1}
((a^*_{ij})_{x_j}(x,u_1)-(a^*_{ij})_{x_j}(x,u_2))
D_iu\,dx\vsp
&=K^{\delta}_1
+K^{\delta}_2
+K^{\delta}_3.
\earr\end{equation}Here, $E_\delta=\{x\in\rr^d;\ |u(x)|\le\delta\}.$ By (H2), it follows that $K^\delta_1\ge0$, while by \eqref{e3.16az}  we have
$$|(a^*_{ij})_u(x,u_1)-(a^*_{ij})_u(x,u_2)|
+|(a^*_{ij})_{x_j}(x,u_1)-(a^*_{ij})_{x_j}(x,u_2)|\le C|u|.$$
Taking into account  that $u_i\in H^1(\rr^d),\ i=1,2,$ and that
\begin{equation}\label{e3.29ac}	 \lim_{\delta\to0}\int_{[|u|\le\delta]} |\nabla u(x)|^2dx=0,\end{equation}we infer that
$\lim\limits_{\delta\to0}K^{\delta}_i=0,$ $i=2,3,$ and so, by \eqref{e2.22a} it follows  that $$\lim_{\delta\to0}\inf I^2_{\delta}\ge0.$$
This yields
$$
|u|_1\le|f|_1,\ \ff\lbb\in(0,\lbb_0).$$

To resume, we have shown so far that under assumptions (H1)--(H3) and (K), for each $f\in L^2$, equation \eqref{e2.3} has, for $\lbb\in (0,\lbb_0)$, a unique solution $u(\lbb,f)\in H^1(\rr^d)$ which satisfies \eqref{e2.16} and \eqref{e2.19}. 

Now, we assume that $a_{ij},b_i$ satisfy (H1)--(H3) only and consider, for $\vp>0$, the functions
\begin{eqnarray}
\quad(a^*_{ij})_\vp(x,u)&\!\!\!=\!\!\!&\int_{\rr}a^*_{ij}(x-\vp y,u-\vp v)\rho(y,v)dv\,dy,\ i,j=1,...,d,\label{e3.28a}\\
b^\vp_{i}(x,u)&\!\!\!=\!\!\!&\int_{\rr}b_{i}(x-\vp y,u-\vp v)\rho(y,v)dv\,dy,\ i,j=1,...,d,\label{e3.31az} 
\end{eqnarray}
where $\rho\in C^\9_0(\rr^d\times\rr)$, $\int_{\rr^d\times\rr}\rho(y,v)dy\,dv=1,$ $\rho\ge0$, is a standard mollifier. Clearly, $(a^*_{ij})_\vp$, $b^\vp_i$ satisfy condition (K). We set $b^\vp=\{b^\vp_i\}^d_{i=1}$. Then, as shown above, the equation
\begin{equation}\label{e3.31aaz} u_\vp-\lbb\dd\sum^d_{i,j=1}D^2_{ij}(a^*_{ij})_\vp(x,u_\vp)+\lbb\ {\rm div}(b^\vp(x,u_\vp)u_\vp)=f\end{equation}
has, for each $\lbb\in(0,\lbb_0)$ and $f\in L^2\cap L^1$, a unique solution $u_\vp=u_\vp(\lbb,f)\in H^1(\rr^d)$ satisfying \eqref{e2.16} and \eqref{e2.19}. Hence
\begin{eqnarray}
\qquad|u_\vp(\lbb,f_1)-u_\vp(\lbb,f_2)|\le|f_1-f_2|_1,\ \ff\vp>0,\ f_1,f_2\in L^2\cap L^1.\label{e3.32z}\\[1mm]
|u_\vp(\lbb,f)|^2_2+\gamma\lbb|\nabla u_\vp(\lbb,f)|^2_2\le C(|f|^2_2+1),\ \ff\vp>0,\ f\in L^2.\label{e3.32az}
\end{eqnarray}
(We note that, by Lemma \ref{l1}, $\lambda_0$ is independent of $\varepsilon$, because
$$\sup\limits_\varepsilon\{|(a^\varepsilon_{ij})_{x_j}|_\infty+|b_i|_\infty;\ i,j=1,...,d\}<\infty.)$$
Now, for $\vp\to0$, it follows by the compactness of $H^1(\rr^d)$ in $L^2_{\rm loc}$ that along a subsequence, again denoted $\vp$, we have
$$u_\vp(\lbb,f)\longrightarrow u\mbox{ strongly in }L^2_{\rm loc}$$ and so, by \eqref{e3.28a}, \eqref{e3.31az}, we have
$$\barr{rcll}
(a^*_{ij})_\vp(x,u_\vp(x))&
\longrightarrow
&a^*_{ij}(x,u(x)),&\mbox{ a.e. }x\in\rr^d,\vsp
b^\vp_i(x,u_\vp(x))& \longrightarrow & b_i(x,u(x)),&\mbox{ a.e. }x\in\rr^d,\earr$$
as $\varepsilon\to0$. Hence, for $\varepsilon\to0$, 
$$\barr{rcll}
D^2_{ij}((a^*_{ij})_\vp(x,u_\vp))&\longrightarrow&D^2_{ij}(a^*_{ij}(x,u))&\mbox{ in }\cald'(\rr^d),\vsp 
{\rm div}(b^\vp(x,u_\vp)u_\vp)&\longrightarrow&{\rm div}(b(x,u)u)&\mbox{ in }\cald'(\rr^d).\earr$$
and so $u=u(\lbb,f)$ is a solution to \eqref{e2.3}. Moreover, by \eqref{e3.32z} it follows that 
\begin{equation}\label{e3.33prim}
|u(\lbb,f_1)-u(\lbb,f_2)|_1\le|f_1-f_2|_1,\ \ff f_1,f_2\in L^2.\end{equation}
Now, we fix $f\in L^1$ and consider a sequence $\{f_n\}\subset L^2$ such that $f_n\to f$ in $L^1$ and consider the corresponding solution $u_n=u(\lbb,f_n)$ to \eqref{e2.3}. 
By \eqref{e3.33prim}, we see that
$$|u_n-u_m|_1\le|f_n-f_m|_1,\ \ff n,m\in\nn.$$Hence, there is $u^*=\lim\limits_{n\to\9}u_n$ in $L^1$. Moreover, by (H1), we see that, for $n\to\infty$,
$$a^*_{ij}(u_n)\to a_{ij}(u^*),\mbox{ a.e. in }\rr^d$$and, since $a_{ij}\in C_b(\rr^d\times\rr)$, we have
$$D^2_{ij}a^*_{ij}(u_n)\to D^2_{ij}a^*_{ij}(u^*)\mbox{ in }\cald'(\rr^d),$$for all  $i,j=1,2,...,d.$ Similarly, by (H3) we see that
$${\rm div}(b(x,u_n)u_n)\to{\rm div}(b(x,u^*)u^*)\mbox{ in }\cald'(\rr^d).$$We have, therefore,
\begin{equation}\label{e3.32azi}
\barr{r}
\dd\sum^d_{i,j=1}D^2_{ij}a^*_{ij}(u_n)-{\rm div}(b(x,u_n)u_n)\!\to\!\sum^d_{i,j=1}D^2_{ij}a^*_{ij}(u^*)
-{\rm div}(b(x,u^*)u^*)\\\mbox{ strongly in }L^1.\earr\hspace*{-2mm}\end{equation}Then, letting $n\to\9$ in equation \eqref{e2.3}, where $f=f_n$, $u=u_n$, we see that $u^*=u(\lbb,f)$ is the solution to \eqref{e2.3}. Moreover, by \eqref{e3.33prim},  the inequality   \eqref{e2.5} follows for all $\lbb\in(0,\lbb_0]$. This means that
$$|(I+\lambda A)^{-1}f_1-(I+\lambda A)^{-1}|_1\le|f_1-f_2|_1,\ \forall f_1,f_2\in L^1,$$for all $\lambda\in[0,\lambda_0).$ By Proposition 3.1 in \cite{1}, this implies that the above inequality holds for all $\lambda>0.$ Hence, \eqref{e2.5} follows for all $\lambda>0$, as claimed.

	As regards \eqref{e2.8},   it first  follows by equation \eqref{e2.3},   where $f\in L^2$ and $u\in H^1(\rr^d)$, by integrating over $\rr^d$. Then, by density, it extends to all of $f\in L^1$. Finally, \eqref{e2.7} for $f\in L^2$, $f\ge0$, follows  by multiplying \eqref{e2.3} with ${\rm sign}(u^-)$ (or, more exactly, by $\calx_\delta(u^-)$ and letting $\delta\to0$) and integrating over $\rr^d$. This completes the proof of Proposition \ref{p1} under hypotheses (H1)--(H3).

	Now, we are ready to formulate the existence theorem for equation \eqref{e1.1}. As mentioned earlier, we shall  represent equation \eqref{e1.1} as the evolution equation \eqref{e1.2} in $\calx=L^1$, where the operator $A$ is defined by \eqref{e2.1}-\eqref{e2.2}. By a {\it weak solution to equation \eqref{e1.1}}, we mean {\it a mild solution to equation \eqref{e1.2},} where $\mathcal{X}=L^1$ and $A$ is the operator defined by \eqref{e2.1}, \eqref{e2.2}.
	
We have
\begin{theorem}\label{t1} Assume that hypotheses {\rm(H1)--(H3)} hold. Then, for each $u_0\in L^1(\rr^d)$, there is a unique weak solution $u=u(\cdot,u_0)\in C([0,\9);L^1)$ to equation \eqref{e1.1}. Moreover, $u$ has the following properties	
\begin{eqnarray}
&|u(t,u^1_0)-u(t,u^2_0)|_1\le|u^1_0-u^2_0|_1,\ \ff u^1_0,u^2_0\in L^1,\ t\ge0,\label{e2.24}\\[2mm]
&u\ge0\mbox{ a.e. in }(0,\9)\times\rr^d\mbox{ if }u_0\ge0\mbox{ a.e. in }\rr^d,\label{e2.26}\\[1mm]
&\dd\int_{\rr^d}u(t,x)dx=\int_{\rr^d}u_0(x)dx, \ff u_0\in L^1,\ t\ge0,\label{e2.27}
\end{eqnarray}
and $u$ is a solution to \eqref{e1.1} in the sense of Schwartz distributions on  \mbox{$ (0,\9){\times}\rr^d,$} $($see~\eqref{1.2}$)$, that~is,
\begin{equation}\label{e2.29a}
\barr{l}
\dd\int^\9_0\int_{\rr^d}(u(t,x)\vf_t(t,x)+\sum^d_{i,j=1} a_{ij}(x,u(t,x))u(t,x)D^2_{ij}\vf(t,x)\\
\qquad+b(x,u)\cdot\nabla_x\vf(t,x)u(t,x))dt\,dx=0,\ \ff\vf\in C^\9_0((0,\9)\times\rr^d).\earr\end{equation}
\end{theorem}

\noindent{\bf Proof.} As mentioned above, the existence of a mild solution $u$ for \eqref{e1.2}, which by our definition is a weak solution to \eqref{e1.1}, follows by the Crandall and Liggett theorem by virtue of   Proposition \ref{p1}, which implies the $m$-accretivity  of the operator $A$ defined by \eqref{e2.1}--\eqref{e2.2}. 
The solution can be  equivalently expressed by the exponential formula
\begin{equation}\label{e2.28}
u(t,u_0)=\lim_{n\to\9}\(I+\frac tn\ A\)^{-n}u_0,\ \ff t\ge0,\ u_0\in\overline{D(A)}=L^1.\end{equation}
Then, by \eqref{e2.5}-\eqref{e2.8}, we get for $u=u(t,u_0)$ the corresponding properties \eqref{e2.24}-\eqref{e2.27} and this completes the proof. In particular, it follows that, if $u_0$ is a   probability density, that is, $u_0dx\in\calp(\rr^d)$, then so is $u(t,u_0)$ for all $t\ge0.$ Note  also that $t\to u(t,u_0)$ is a continuous semigroup of nonexpansive operators in the space~$L^1$. As regards \eqref{e2.29a}, it follows by letting $h\to0$ in the equation
$$\begin{array}{r}
\displaystyle\int^\infty_0\int_{\rr^d}
\Big(u_h(t,x)(\varphi(t+h,x)-\varphi(t,x))
+\sum^d_{i,j=1}a_{ij}(x,u_h(t,x))u_h(t,x)\vsp 
D^2_{ij}\varphi(t,x)+b(x,u_h(t,x))\cdot\nabla_x\varphi(t,x)u_h(t,x)\Big)dt\,dx=0,\vsp \hfill\forall\varphi\in C^\infty_0((0,\infty)\times\rr^d).\end{array}$$

\begin{remark}\label{r3.4}\rm Assumptions $a_{ij}\in C^2(\rr^d\times\rr)$ and $b_i\in C^1(\rr^d\times\rr)$ in (H1)--(H3) were necessary for the density of $D(A)$ in $L^1$. Otherwise, it suffices to take only $a_{ij}\in C^1(\rr^d\times\rr)$. In the special case $a_{ij}\equiv\beta\delta_{ij}$, the density of $D(A)$ follows, however, under the weaker condition $\beta\in C^1(\rr^d)$  (see \cite{6ab}).\end{remark}
	
\subsection{Existence for degenerate FPEs}

 We consider here the equation
\begin{equation}\label{e3.1}
\barr{l}
u_t-\dd\sum^d_{i,j=1} D^2_{ij}(a_{ij}(u)u)+\dd\sum^d_{i=1}D_i(b_i(u)u)=0\mbox{ in }\cald'((0,\9)\times\rr^d),\vsp
u(0,x)=u_0(x), x\in\rr^d,\earr\end{equation}
where $a_{ij}$ and $b_i$ satisfy hypotheses (H1)$'$--(H3)$'$.

Consider the operator $A_1:D(A_1)\subset L^1\to L^1$ defined by
\begin{equation}\label{e3.4}
\barr{l}
A_1u=-\dd\sum^d_{i,j=1}D^2_{ij}(a_{ij}(u)u)
+\sum^d_{i=1}D_i(b_i(u)u)\mbox{ in }\cald'(\rr^d),\vsp
D(A_1)=\left\{u\in L^1;-\dd\sum^d_{i,j=1}D^2_{ij}(a_{ij}(u)u)+
\sum^d_{i=1}D_i(b_i(u)u)\in L^1\right\}.\earr\end{equation}
We have

\begin{lemma}\label{l4} Assume that {\rm(H1)$'$--(H3)$'$} hold. Then the operator $A_1$ is $m$-accretive in $L^1$.
\end{lemma}

\n{\bf Proof.} One should prove that, for each  $\lbb\in(0,\lbb_0)$     and   $f\in L^1$, the equation
\begin{equation}\label{e3.5}
u-\lbb\sum^d_{i,j=1}D^2_{ij}(a_{ij}(u)u)
+\lbb\sum^d_{i=1}D_i(b_i(u)u)=f\mbox{ in }\cald'(\rr^d)\end{equation}has a unique solution $u=u(\lbb,f)$ which satisfies the estimate
\begin{equation}\label{e3.6a}
|u(\lbb,f_1)-u(\lbb,f_2)|_1\le |f_1-f_2|_1,\ \ff f_1,f_2\in L^1.\end{equation}

We set, for each $\vp>0$,
\begin{equation}
a^\vp_{ij}(r)=a_{ij}(r)+\vp\delta_{ij},\   i,j=1,...,d,\ r\in\rr, \label{e3.4a}
\end{equation}
where $\delta_{ij}$ is the Kronecker symbol.
Then, we approximate \eqref{e3.5} by
\begin{equation}
\label{e3.4aaa}
u-\lbb\sum^d_{i,j=1}D^2_{ij}(a^\vp_{ij}(u)u)
+\lbb\sum^d_{i=1}D_i(b_i(u)u)=f\mbox{ in }\cald'(\rr^d).
\end{equation}
Equivalently,
\begin{equation}\label{e3.42prim}
u+\lbb A^\vp_1(u)=f,\end{equation}where
$$\barr{rcl}
A^\vp_1(u)&=&-\dd\sum^d_{i,j=1}D^2_{ij}(a^\vp_{ij}(u)u)+
\sum^d_{i=1}D_i(b_i(u)u),\ \ff u\in D(A^\vp_1),\vsp
D(A^\vp_1)&=&\left\{u\in L^1;\ -\dd\sum^d_{i,j=1}D^2_{ij}(a^\vp_{ij}(u)u)
+\dd\sum^d_{i=1}D_i(b_i(u)u)\in L^1\right\}.\earr$$
We shall prove   that, for each $f\in L^1$, there is a solution $u=u_\vp(\lbb,f)$ satisfying \eqref{e3.6a} for $0<\lbb<\lbb_0.$

Since $a^\vp_{ij}$ and $b_i$ satisfy, for each $\vp>0$, hypotheses (H1)--(H3),   Proposition \ref{p1}  implies the existence  of a solution $u_\vp=u_\vp(\lbb,f)$ to \eqref{e3.4aaa} in $L^1(\rr^d)$ for each $f\in L^2$ if $0<\lbb\le\lbb^\vp_0=\frac C\vp,$ with $C$ independent of $\vp$.

Moreover, one has \begin{equation}\label{e3.42b}
|u_\vp(\lbb,f_1)-u_\vp(\lbb,f_2)|_1\le|f_1-f_2|_1,\ \ff f_1,f_2\in L^2,\ \lbb\in(0,\lbb^\vp_0).\end{equation}
Then,  by   density, $u_\vp(\lbb,f)$ extends as solution to \eqref{e3.4aaa} for all $f\in L^1$.

Note also that, by   \eqref{e2.5}--\eqref{e2.8}, we have, for all $\vp>0$ and $\lbb\in(0,\lbb^\vp_0)$,  \begin{eqnarray}
&\dd\int_{\rr^d}(I+\lbb A^\vp_1)\1 f\,dx=\int_{\rr^d}f\,dx,\ \ \ff f\in L^\9,\label{e3.7}\\[2mm]
&(I+\lbb A^\vp_1)\1f\ge0,\mbox{ a.e. in }\rr^d\mbox{ if }f\ge0,\mbox{ a.e. in }\rr^d,\label{e3.8}\end{eqnarray}while  \eqref{e3.42b} yields
\begin{equation}
 |I{+}\lbb A^\vp_1)\1f_1{-}(I{+}\lbb A^\vp_1)\1f_2|_1\le|f_1{-}f_2|_1,\, \ff f_1,f_2\in L^1,\, \vp>0.\label{e3.7a}
\end{equation}

Though \eqref{e3.7}--\eqref{e3.7a} were proved only  for $0<\lbb\le\lbb^\vp_0$, it can be shown, however, as mentioned earlier, that $(I+\lbb A^\vp_1)\1$ extends to all  $\lbb>0$ by a well known argument based on the resolvent equation
$$(I+\lbb A^\vp_1)\1 f=(I+\lbb_0A^\vp_1)\1\(\frac{\lbb^\vp_0}\lbb\,f+\(1-\frac{\lbb^\vp_0}\lbb\)(I+\lbb A^\vp_1)\1f\),\ \lbb>\lbb^\vp_0.$$(See \cite{1}, Proposition 3.3.)

Now, we are going to let $\vp\to0$ in \eqref{e3.4aaa}. We set, for $f\in L^1$ and  the solution $u_\vp$ to \eqref{e3.4aaa},
$$u^\vp_h(x)=u_\vp(x+h)-u_\vp(x),\ f_h(x)=f(x+h)-f(x),\ x,h\in\rr^d.$$
Since $a^\vp_{ij}$ and $b^\vp_i$ are independent of $x$, we see that $x\to u^\vp(x+h)$ is the solution to \eqref{e3.4aaa} for $f(x)=f(x+h)$. Then, by \eqref{e3.7a}, it follows that
$$|u^\vp_h|_1\le |f_h|_1,\ \ff h\in\rr^d,\ \vp>0.$$
By the Kolmogorov compactness theorem (see, e.g., \cite{6a}, p.~111), it follows that $\{u^\vp\}$ is compact in $L^1_{\rm loc}(\rr^d)$ and so, along a subsequence, $$\mbox{$u_\vp\to u$ strongly in $L^1_{\rm loc}(\rr^d)$ for $\vp\to0$.}$$Since $|u_\vp|_1\le C,\ \ff\vp>0,$ it follows via Fatou's  lemma that $u\in L^1$. Letting $\vp\to0$ in \eqref{e3.4aaa},   where $u=u_\vp$, and taking into account that
$$a^\vp_{ij}(u_\vp)u_\vp\to a_{ij}(u)u,\ \ b_i(u_\vp)u_\vp\to b_i(u)u,\mbox{ a.e. in }\rr^d,$$
while by (H1)$'$, (H3)$'$,
$$|a^\vp_{ij}(u_\vp)|+|b_{ij}(u_\vp)|\le C,\mbox{ a.e. in }\rr^d,$$where $C$ is independent of $\vp$, we see that $u$ is a solution to \eqref{e3.4aaa} and so $u=(I+\lbb A_1)\1f.$ Moreover, letting $\vp\to0$ in \eqref{e3.7}--\eqref{e3.7a}, we see that
\begin{eqnarray}
&\ |(I{+}\lbb A_1)\1f_1{-}(I{+}\lbb A_1)\1f_2|_1{\le}|f_1{-}f_2|_1,\, \ff\lbb>0,\, f_1,f_2\in L^1,\label{e3.10az}\\[2mm]
&\dd\int_{\rr^d}(I+\lbb A_1)\1f\,dx=\int_{\rr^d}f\,dx,\ \ff f\in L^1,\ \lbb>0,\label{e3.10aaz}\\[2mm]
&(I+\lbb A_1)\1f\ge0,\mbox{ a.e. in }\rr^d\mbox{ if }f\ge0,\mbox{ a.e. in }\rr^d.\label{e3.10aaaz}
\end{eqnarray}
Then, by the Crandall and Liggett existence theorem, for each $u_0{\in}\overline{D(A_1)}{=}L^1$, the differential equation
\begin{equation}\label{e3.9}
\barr{l}
\dd\frac{du}{dt}+A_1u=0,\ t>0,\vsp
u(0)=u_0,\earr\end{equation}{\it has a unique mild solution $u\in C([0,\9);L^1)$ in the sense of \eqref{e1.4}--\eqref{e1.6}.}

As in the previous case, this mild solution is, in fact, a solution to the Fokker-Planck equation \eqref{e3.1} in the sense of Schwartz distributions (cf.~\eqref{e2.29a}). 
We have, therefore, the following existence result.

\begin{theorem}\label{t2} Under hypotheses {\rm(H1)$'$--(H3)$'$}, for  each $u_0\in L^1$, there is a unique weak solution $u=u(t,u_0)\in C([0,\9);L^1)$ to equation \eqref{e3.1}. Moreover, this solution satisfies \eqref{e2.24}--\eqref{e2.27} and is a solution to \eqref{e3.1} in the sense of Schwartz distributions on $(0,\9)\times\rr^d$, i.e., in the sense of \eqref{e2.29a} or \eqref{1.2}.
\end{theorem}

\begin{remark}\label{r1}\rm In particular, Theorems \ref{t1} and \ref{t2} imply the existence of   a solution $u$ in the sense of Schwartz distributions on $(0,\9)\times\rr^d$ for equation \eqref{e1.1}. Moreover, $u:[0,\9)\to L^1$ is continuous. In some special cases, these two properties  are sufficient to characterize such solutions to \eqref{e1.1}. In fact, this is the case if (see \cite{6}) $b\equiv0$ and
		$$a_{ij}(x,u)u=\delta_{ij}\beta(u)u,\ \ff u\in\rr,\ i,j=1,...,d,$$ where  $\beta$ is a continuous monotonically nondecreasing function because, in this case, one has the uniqueness of distributional solutions  \mbox{$u{\in} L^\9((0,\9){\times}\rr^d)$} $\cap C([0,\9);L^1)$. Such a result remains, however, open for   general  Fokker-Planck equations as in \eqref{e1.1}.	\end{remark}

\begin{remark}\label{r3} \rm In the special case $a_{ij}=\delta_{ij}$, the weak solution $u$ given by Theorem \ref{t2} is an  entropic solution in sense of S. Kruzkov \cite{19a} for equation \eqref{e1.1}.  In the present case, the solution $u$ given by Theorem \ref{t2} is a {\it mild} solution to \eqref{e1.1} defined, as in the previous case, by the finite difference scheme \eqref{e1.4}--\eqref{e1.6}. It is, of course, a continuous in $t$ distributional solution to \eqref{e1.1}, but we do not know if it is unique within this class. In fact, we should mention that the solutions $u$ given by Theorems \ref{t1} and \ref{t2} are unique in the class of mild solutions generated by the operator $A$ and not in the class of distributional or entropic solutions in the sense of Kruzkov.\end{remark}

\section{Solution of the McKean-Vlasov SDE}\label{s4}
\setcounter{equation}{0}

Consider the following McKean-Vlasov SDE for $T\in(0,\9)$
\begin{equation}\label{e4.1}
\barr{rcl}
dX(t)&=&b\(X(t),\dd\frac{d\call_{X(t)}}{dx}\,(X(t))\)dt\vsp&&
+\sqrt{2}\ \sigma\(X(t),\dd\frac{d\call_{X(t)}}{dx}\,(X(t))\)dW(t),\ 0\le t\le T,\vsp
X(0)&=&\xi_0,\earr \end{equation}on $\rr^d$, where $W(t),\ t\ge0,$ is an $(\calf_t)_{t\ge0}$-Brownian motion on a probability space $(\ooo,\calf,P)$ with normal filtration $(\calf_t)_{t\ge0}$ and $\xi_0:\ooo\to\rr^d$ is $\calf_0$-measurable such that $$P\circ\xi^{-1}_0(dx)=u_0(x)dx.$$ Furthermore, $b=(b_1,...,b_d):\rr^d\times\rr\to\rr^d$ and $\sigma:\rr^d\times\rr\to L(\rr^d;\rr^d)$ are measurable.

Let $a_{ij}:=2(\sigma\sigma^T)_{ij},\ 1\le i,j\le d.$ Then, as an immediate consequence of Section \ref{s2} and Theorems \ref{t1} and \ref{t2}, respectively, we obtain the following.

\begin{theorem}\label{t4.1} Suppose that $a_{ij},b_i,$ $1\le i,j\le d,$ satisfy either {\rm(H1)--(H3)} or  {\rm(H1)$'$--(H3)$'$}. Then there exists a $($in the probabilistic sense$)$  weak solution to DDSDE \eqref{e4.1}. Furthermore, for the solution $u$ in Theorem {\rm\ref{t1}} and {\rm\ref{t2}}, respectively,    with $u(0,\cdot)=u_0$, we have the "probabilistic representation"
	$$u(t,x)dx=P\circ X(t)\1(dx),\ t\ge0.$$
\end{theorem}

\begin{remark}\label{r4.1}\rm \vspace*{-2mm}\
	\begin{itemize}\item[(i)] In the   case  where in \eqref{e4.1} we have $a_{ij}(x,u)=\delta_{ij}\beta(u),$ $ 1\le i,j\le d$, and $\beta:\rr\to2^\rr$ is maximal monotone  with  $\sup\{|s|:s\in\beta(r)\}\le C|r|^m,$ $r\in\rr$, for some $C,m\in[0,\9)$ and $b$ satisfies (H3)$'$, then the above theorem was already proved in \cite{3}. The special case where, in addition, $b\equiv0$, $d=1$ and $m=4$, was proved in \cite{5} if $\beta(r)/r$ is nondegenerate at $r=0$ and in \cite{13*} including the degenerate case.\vspace*{-2mm}
		\item[(ii)] The special case $d=1$, $b\equiv0$, $a_{ij}(x,u)=\delta_{ij}\beta(u),$ $1\le i,j\le d$, with $\beta(r):=r|r|^{m-1},$ $r\in\rr$, for some $m\in(1,\9)$, was proved in \cite{16*}.\vspace*{-2mm}
		\item[(iii)] \cite{5prim} contains an analogous result as in \cite{5}, \cite{13*} in the case where a linear multiplicative noise is added to the nonlinear FPE, which thus becomes a stochastic porous media equation. \end{itemize}\end{remark}
	
	Our final remark concerns the uniqueness of the time marginal of solutions to \eqref{e4.1}.
	
\begin{remark}\label{r4.2}\rm If $b\equiv0$ and $a_{ij}(x,u)=\delta_{ij}\beta(u),$ $1\le i,j\le d$, and $\beta:\rr\to\rr$ is continuous, nondecreasing  and $\beta(0)=0$, then \eqref{e1.1} has a unique solution among all the solutions in $(L^\9\cap L^1)((0,T)\times\rr^d)$ by the main result in \cite{6}. Hence, obviously, we have uniqueness of the time marginals for weak solutions to \eqref{e4.1} among all the solutions of \eqref{e4.1} whose time marginals have densities in $L^\9((0,T)\times\rr^d)$.\end{remark}

\n{\bf Acknowledgements.} This work    was supported by the DFG through CRC  1283. V. Barbu was partially supported by the   CNCS-UEFISCDI project PN-III-P4-ID-PCE-2015-0011.


\begin{thebibliography}{nn}
	
	\bibitem{1}  Barbu, V., {\it Nonlinear Differential Equations of Monotone Type in Banach Spaces}, Springer, 2010.
	
	\bibitem{2}  Barbu, V., Generalized solutions to nonlinear Fokker-Planck equations, {\it J.~Diff. Equations}, 261 (2016), 2446-2471.
	
\bibitem{2a}  Barbu, V., R\"ockner, M., Nonlinear   Fokker-Planck equations driven by Gaussian linear multiplicative noise, {\it J.~Differential Equations} (2018).

\bibitem{3}  Barbu, V., R\"ockner, M., Probabilistic representation for solutions to nonlinear Fokker-Planck equation, {\it SIAM J. Math. Anal.}, 50 (2018), 2588-2607.
    
\bibitem{4}  Barbu, V., R\"ockner, M., Russo, F., Probabilistic representation for solutions of an irregular porous media type equation: the degenerate case, {\it Probab. Theory Rel. Fields,} 15 (2011), 1-43.


\bibitem{5prim} Barbu, V., R\"ockner, M., Russo, F., Doubly probabilistic representation for the stochastic porous media equation, {\it Ann. Inst. Henri Poinfcar\'e}, 53 (4) (2017), 1061-1090.

\bibitem{6ab} Barbu, V., R\"ockner, M., The evolution to equilibrium of solutions to nonlinear Fokker-Planck equations, arXiv:1904.08291.



\bibitem{16*}
Benachour, S., Chassaing, P., Roynette, B., Vallois, P.,
Processus associ\'es \`a l'\'equation des milieux poreux,
{\it Ann. Scuola Norm. Sup. Pisa Cl. Sci.}, (4) (1997), 23 (4) (1996), 793--832.


\bibitem{5}
Blanchard, Ph., R\"ockner, M., Russo, F., Probabilistic representation for solutions of an irregular porous media equation, {\it Ann. Probab.}, 38 (2010), 1870-1900.

\bibitem{13*}
Bogachev, V.I., Krylov, N.V., R\"ockner, M.,  Stanislav~V.
Shaposhnikov, S.V.,
{\it Fokker-{P}lanck-{K}olmogorov equations},
Mathematical Surveys and Monographs, vol. 207,
American Mathematical Society, Providence, RI, 2015.

    \bibitem{6a} Brezis, H., {\it Functional Analysis, Sobolev Spaces and Partial Differential Equations}, Springer, New York. Dordrecht. Heidelberg. London, 2011.

\bibitem{6} Brezis, H., Crandall, M.G., Uniqueness of solutions of the initial-value problem for $u_t-\Delta\beta(u)=0,$ {\it J.~Math. Pures et Appl.}, 58 (1979), 153-163.

\bibitem{9*}
dos Reis, G., Smith, G., Tankov, P.,
{I}mportance sampling for {M}c{K}ean-{V}lasov {S}{D}{E}s,
arXiv:1803.09320, 2018.

\bibitem{14*}
Figalli, A.,
Existence and uniqueness of martingale solutions for {SDE}s with
rough or degenerate coefficients,
{\it J. Funct. Anal.}, 254 (1) (2008), 109--153.


\bibitem{8a} Franck, T.D., {\it Nonlinear Fokker-Planck Equations. Fundamentals and Applications}, Springer, Berlin. Heidelberg. New York, 2005.


\bibitem{4*}
Funaki, T.,
A certain class of diffusion processes associated with nonlinear
parabolic equations,
{\it Z. Wahrsch. Verw. Gebiete}, 67 (3) (1984), 331-348.


\bibitem{8*}
Hammersley, W., Šiška, D., Szpruch, L.,
{M}c{K}ean-{V}lasov {S}{D}{E}s under measure mependent Lyapunov
conditions,
arXiv:1802.03974, 2018.

\bibitem{11*}
Huang, X., R\"ockner, M., Wang, F.-Y.,
{N}onlinear {F}okker-{P}lanck equations for probability measures
on path space and path-distribution dependent {S}{D}{E}s,
arXiv:1709.00556, 2017.

\bibitem{10*}
Huang, X.,  Wang, F.-Y.,
{D}istribution dependent {S}{D}{E}s with singular coefficients,
arXiv:1805.01682, 2018.


\bibitem{19a} Kruzkov, S., First order quasilinear equations with several independent variables, {\it Sbornic: Mathematics}, 10 (2) (1970), 217-243.

\bibitem{1*}
McKean, Jr., H.P.,
A class of {M}arkov processes associated with nonlinear parabolic
equations,
{\it Proc. Nat. Acad. Sci. U.S.A.}, 56 (1966), 1907-1911.

\bibitem{2*}
McKean, Jr., H.P.,
Propagation of chaos for a class of non-linear parabolic equations.
In {\it Stochastic {D}ifferential {E}quations} ({L}ecture {S}eries in
{D}ifferential {E}quations,  {S}ession 7, {C}atholic {U}niv., 1967), pages
41-57. Air Force Office Sci. Res., Arlington, Va., 1967.


\bibitem{12*}
Mehri, S., Scheutzow, M., Stannat, W.,  Zangeneh, B.Z.,
{P}ropagation of {C}haos for {S}tochastic {S}patially {S}tructured
{N}euronal {N}etworks with {F}ully {P}ath {D}ependent {D}elays and {M}onotone
{C}oefficients driven by {J}ump {D}iffusion {N}oise,
arXiv:1805.01654, 2018.


\bibitem{7*}
Mishura, Yu.S.,     Veretennikov, A.Yu.,
{E}xistence and uniqueness theorems for solutions of
{M}c{K}ean--{V}lasov stochastic equations,
arXiv:1603.02212, 2016.

\bibitem{13*a}
R\"ockner, M., Russo, F., Uniqueness of a class of stochastic Fokker-Planck and porous media equations, {\it J. Evol. Equ.}, 17 (3) (2017), 1049-1062.

\bibitem{5*}
Scheutzow, M.,
Uniqueness and nonuniqueness of solutions of {V}lasov-{M}c{K}ean
equations,
{\it J. Austral. Math. Soc. Ser. A}, 43 (2) (1987), 246-256.


\bibitem{15*}
Stroock D.W., Srinivasa Varadhan, S.R.,
{\it Multidimensional diffusion processes}, 
Springer, Berlin. New York, 1997.

\bibitem{3*}
Sznitman, A.-S.,
Nonlinear reflecting diffusion process, and the propagation of chaos
and fluctuations associated,
{\it J. Funct. Anal.}, 56 (3) (1984), 311-336.

\bibitem{8} Trevisan, D., Well-posedness of multidimensional diffusion processes, {\it Electronic J. Probab.}, 21 (2016), 1-41.

\bibitem{6*}
Wang, F.-Y.,
Distribution dependent {SDE}s for {L}andau type equations,
{\it Stochastic Process. Appl.}, 128 (2) (2018), 595-621.


\end{thebibliography}
\end{document}